\documentclass[12pt]{article}
\usepackage{amsmath}
\usepackage{amssymb}
\usepackage{amsfonts}
\newtheorem{thm}{Theorem}[section]
\newtheorem{pro}[thm]{Proposition}

\newtheorem{dfn}[thm]{Definition}
\newtheorem{fct}[thm]{Fact}

\newtheorem{rmk}[thm]{Remark}

=22cm\headheight=0cm\topskip=0cm\topmargin=0cm

\def\qed{\hfill$\Box$}

\newcommand{\M}{\sf M}

\begin{document}
\def\dis{\displaystyle}

\begin{center}
{\LARGE{ NIP formulas and Baire 1 definability }} \vspace{10mm}

{\large{ Karim Khanaki}} \vspace{3mm}

Department of science,\\
  Arak University of Technology,\\
P.O. Box 38135-1177, Arak, Iran; \\
Email: k.khanaki@gmail.com

\end{center}

\begin{abstract}
In this short note, using results of Bourgain, Fremlin, and
Talagrand  \cite{BFT}, we show that for a countable structure
$M$, a saturated elementary extension $M^*$ of $M$ and a formula
$\phi(x,y)$ the following are equivalent:
\begin{itemize}
             \item [(i)]  $\phi(x,y)$ is NIP on $M$ (in the sense of Definition~\ref{NIP-formula}).

             \item [(ii)]  Whenever $p(x)\in S_\phi(M^*)$ is finitely satisfiable in $M$ then it is Baire 1 definable over $M$ (in sense of Definition~\ref{Baire 1 definability}).
\end{itemize}
\end{abstract}





\section{Introduction} \label{}
A note of Anand Pillay \cite{P} encourages us to write the present
note. In his work, Pillay shows that a formula $\phi(x,y)$ on a
structure $M$ has NOT order property iff very type $p$ in
$S_\phi(M)$  has an extension to a type $p'$ in $S_\phi(M^*)$
(where $M^*$ is a saturated elementary extension of $M$) such that
$p'$ is both definable over and finitely satisfiable in $M$. He
pointed out that a model theoretic proof of this equivalence was
given by himself in \cite{P82}. In fact this is a well known
result essentially due to A. Grothendieck \cite{Gro} which
asserts that $\phi(x,y)$ has NOT order property on $M$ iff the set
$\{\phi(a,y):S_{\tilde{\phi}}(M)\to\{0,1\}:a\in M \}$ is
relatively weakly compact in $C(S_{\tilde\phi}(M))$, the set of
all continuous  0-1  valued functions on the space of
$\phi$-types. Moreover, he presented a model-theoretic proof of
Grothendieck's result in \cite{P}. Of course, this connection was
formerly known by some people, including Ben Yaacov
\cite{Ben-Gro}.

On the other hand, recently some authors \cite{Iba, K, S2}
observed the connection between non independence property (NIP)
and a result of Bourgain, Fremlin, and Talagrand in \cite{BFT}.
Roughly speaking, like stability, a formula $\phi$ is NIP iff a
set of continuous functions has a `good' topological property,
namely relative sequential compactness.

 In this short note, we present a result for NIP formulas similar to
 the result of Pillay \cite{P} for NOP formulas. Our approach is
 `double local', i.e., we do not study an NIP theory $T$ or  an NIP formula
 $\phi$, but we study a formula $\phi$ which is NIP on a structure
 $M$. With these assumptions, we are in the realm of mathematical analysis.
Our main observation is just a translation of the main result in
\cite{BFT} into the language of model theory. For two reasons we
do not provide a model-theoretic proof: (1) the result in
\cite{BFT} is very general, i.e., it provides a similar result for
$[0,1]$ valued model theory, namely continuous logic \cite{BBHU}.
Also, combinatorial aspects of the proof in $[0,1]$ version are
much more difficult than 0-1 valued logic. (2) This is an
invitation for model theorists who pay more attention to the
relationship between model theory and Functional Analysis.

\section{Definitions and proof}
\begin{dfn} \label{NIP-formula}
{\em  Let $\M$ be a structure, and $\phi(x,y)$ a formula.
 We say that $\phi(x,y)$ is NIP on $M\times M$ (or on $M$) if
for each sequence $(a_n)\subseteq M$,  there are some {\em
finite} disjoint subsets $E,F$ of ${\Bbb N}$ such that
$${\M}\nvDash\exists y\big( \bigwedge_{n\in
 E}\phi(a_n,y)~\wedge~\bigwedge_{n\in F}\neg\phi(a_n,y)\big).$$}
\end{dfn}

There are more equivalent notions of NIP in the frame of
continuous logic (see \cite{K-Banach}). This notion was used by
some authors in \cite{Iba} and \cite{K}. (Note that since we want
to study countable structures, in Definition~\ref{NIP-formula},
the subsets $E,F$ must be finite.)

\begin{fct} \label{NIP=RSC}  Let $\M$ be a structure, and $\phi(x,y)$ a formula.
The the following are equivalent.
\begin{itemize}
             \item [{\em (i)}] $\phi(x,y)$ is NIP on $M$.
             \item [{\em (ii)}] For each sequence $\phi(a_n,y)$ in the set
                        $A=\{\phi(a,y):S_y(M)\to\{0,1\}~|~a\in M\}$, where $S_y(M)$ is
                        the space of all complete types on $M$ in the variable $y$,
                        there are some {\em arbitrary} disjoint subsets $E,F$ of ${\Bbb N}$ such that
                         $$\Big\{y\in S_y(M):\big( \bigwedge_{n\in E}\phi(a_n,y)=0\big)\wedge\big(\bigwedge_{n\in F}\phi(a_n,y)=1\big)\Big\}=\emptyset.$$
Where $\phi(a,b)=1$ if $\vDash\phi(a,b)$ and $\phi(a,b)=0$ if
$\vDash\neg\phi(a,b)$.
             \item [{\em (iii)}] Every sequence $\phi(a_n,y)$ in $A$ has a convergent subsequence.
\end{itemize}
\end{fct}

{\bf Proof.}  By the compactness of the type space,
(i)~$\Leftrightarrow$~(ii). The equivalence
(ii)~$\Leftrightarrow$~(iii) is a well known result due to
Rosenthal (see Rosenthal's lemma in \cite{K-Banach}).~\qed

\medskip
Now, we present main results of Bourgain, Fremlin, and Talagrand
\cite{BFT} which will be used in this note.

\begin{fct}[\cite{BFT}, Corollary~4G]  \label{Polish-compact}
Let $X$ be a Polish space, and $A\subseteq C(X)$ pointwise
bounded set. Then the following are equivalent:
\begin{itemize}
             \item [{\em (i)}] $A$ is relatively compact in the set of all Baire 1 functions on $X$, denoted by $B_1(X)$.
              \item [{\em (ii)}] $A$ is relatively sequentially compact in
              ${\mathbb{R}}^X$, i.e., every sequence in $A$ has
              a pointwise convergent subsequence in ${\Bbb R}^X$.
\end{itemize}
\end{fct}

A regular Hausdorff space X is {\em angelic} if (i) every
relatively countably compact set, i.e., a set so that every
sequence of it has a cluster point in $X$, is relatively compact,
(ii) the closure of a relatively compact set is precisely the set
of limits of its sequences. The following is the principal result
of Bourgain, Fremlin and Talagrand.

\begin{fct}[\cite{BFT}, Theorem~3F] \label{Polish-angelic}
If $X$ is a Polish space, then $B_1(X)$ is angelic under the
topology of pointwise convergence.
\end{fct}

Recall that a function $\psi$ on a topological space is Baire 1 if
it is the pointwise limit of a sequence of continuous functions
on $X$.

\begin{dfn}[Baire 1 definability] \label{Baire 1 definability}  Let $M$ be a countable structure,
 $M^*$ an elementary extension of it, and $\phi(x,y)$ a
formula. A type $q\in S_\phi(M^*)$ is Baire 1 definable over $M$
if there is a sequence $\phi(a_n,y), a_n\in M$, such that, for
all $b\in M^*$, $q\vdash\phi(x,b)$ iff eventually
$\vDash\phi(a_n,b)$, i.e., for each $b\in M^*$, there is some $k$
such that for all $n\geqslant k$, $\vDash\phi(a_n,b)$.
\end{dfn}

Equivalently, if $\phi(a_n,y):S_{\tilde\phi}(M^*)\to\{0,1\}$ and
let $\psi(y):=\lim_n\phi(a_n,y)\upharpoonright_{M^*}$, then for
all $b\in M^*$, $q\vdash\phi(x,b)$ iff $\vDash\psi(b)$.

Recall that for a countable structure $M$, the space $S_\phi(M)$
is a Polish space, since the set $M_0=\{tp_\phi(a/M):a\in M\}$ is
dense in $S_\phi(M)$.

Similar to \cite{P}, we let
${\tilde\phi}(x,y)=\phi^*(x,y)=\phi(y,x)$.

As previously mentioned in introduction, the following
observation is just a translation of the main result in
\cite{BFT} into the language of model theory.

\begin{pro} \label{NIP=Borel}
Let $M$ be a countable structure, $M^*$ a saturated elementary
extension and $\phi(x,y)$ a formula. Then the following are
equivalent:
\begin{itemize}
             \item [{\em (i)}] $\phi(x,y)$ is NIP on $M$.

 \item [{\em (ii)}] Whenever $p(x)\in S_\phi(M^*)$ is finitely satisfiable in $M$ then it is Baire 1 definable over
 $M$.
\end{itemize}
\end{pro}

{\bf Proof.}  (i) $\Rightarrow$ (ii): This is equivalent to
saying  that for every net $(a_i)_{i\in I}$ of elements of $M$,
if $tp_\phi(a_i/M^*)\to p$ in the logic topology, then the type
$p$ is Baire 1 definable over $M$. Indeed, let
$A=\{\phi(a,y):S_{\tilde\phi}(M)\to\{0,1\}:a\in M\}$. Set
$\psi(y)=\lim_i\phi(a_i,y)$ for all $y\in S_{\tilde\phi}(M)$.
Recall that $S_{\tilde\phi}(M)$ is Polish.  By
Facts~\ref{NIP=RSC}, \ref{Polish-compact} and
\ref{Polish-angelic}, there is a sequence $(a_n)\subseteq M$ such
that the sequence $\phi(a_n,y):S_{\tilde\phi}(M)\to\{0,1\}$
pointwise converges to $\psi$. (Note that $\psi$ is Baire 1, but
may not be continuous.)  Now, since $p$ is $M$-invariant, for all
$b\in M^*$, we have
 $\phi(p,b)=\phi(p,tp_{\tilde\phi}(b/M))=\psi(tp_{\tilde\phi}(b/M))=\lim_n\phi(a_n,tp_{\tilde\phi}(b/M))$.
So,
$\lim_n\phi^{M^*}(a_n,b)=\lim_n\phi(a_n,tp_{\tilde\phi}(b/M))=\phi(p,b)$.
To summarize, for all $b\in M^*$,
 $p\vdash\phi(x,b)$ iff $\lim_n\phi^{M^*}(a_n,b)=\phi(p,b)$.
On the other word, $p$ is Baire 1 definable over $M$, by the
Baire 1 function $\psi^{M^*}(y)=\lim_n\phi^{M^*}(a_n,y)$ for all
$y\in M^*$.

(ii) $\Rightarrow$ (i): Indeed, we show that the set $A$ is
relatively  compact in $B_1(S_{\tilde\phi}(M))$. Suppose that
$\phi(a_i,y)\to\psi(y)$ where $a_i$ is a net in $M$. Since $M^*$
is a saturated elementary extension of $M$, all types in
$S_{\tilde\phi}(M)$ are realised in $M^*$. Suppose that
$tp_\phi(a_i/M^*)\to p$ in the logic topology. Clearly, $p\in
S_\phi(M^*)$ and by (ii), it is  Baire 1 definable. Therefore,
there is a sequence $a_n\in M$ such that
$\lim_n\phi(a_n,y)=\phi(p,y)$ for all $y\in M^*$. Therefore,
$\psi(y)=\lim_n\phi(a_n,y)$ for all $y\in S_{\tilde\phi}(M)$. So,
$\psi:S_{\tilde\phi}(M)\to\{0,1\}$ is Baire 1, and hence $\bar
A\subseteq B_1(S_{\tilde\phi}(M))$. So, $A$ is relatively compact
in $B_1(S_{\tilde\phi}(M))$, and therefore by Facts~\ref{NIP=RSC}
and \ref{Polish-compact}, $\phi$ is NIP.~\qed

\begin{rmk} Note that by some straightforward adaptions, all statements of the present note hold
in the framework of continuous logic \cite{BBHU}. One of the
reasons that we did not provide model-theoretic proofs was the
possibility of easy generalizations and adaptations of these
statements for continuous logic using the language of Analysis.
\end{rmk}


\bigskip\noindent
 {\bf Acknowledgements.} I want to thank Anand Pillay for
his helpful comments and his detailed corrections.

 I would like to thank the Institute for Basic
Sciences (IPM), Tehran, Iran. Research partially supported by IPM
grants nos 93030032 and 93030059.


\end{document}